\documentclass[12pt]{article}

\usepackage{amsfonts}			
\newcommand{\RR}{{\mathbb R}}
\newcommand{\CC}{{\mathbb C}}
\newcommand{\HH}{{\mathbb H}}
\newcommand{\OO}{{\mathbb O}}
\newcommand{\KK}{{\mathbb K}}
\newcommand{\TT}{{\mathbb T}}
\newcommand{\FF}{{\mathbb F}}
\DeclareMathSymbol\blacksquare{\mathord}{AMSa}{"04}	

\newcommand{\Aa}{a}			
\newcommand{\Ab}{b}
\newcommand{\Ac}{c}
\newcommand{\Ap}{d}
\newcommand{\Am}{e}
\newcommand{\An}{f}

\newcommand{\Vx}{x}
\newcommand{\Vy}{y}
\newcommand{\Vz}{z}
\newcommand{\Ox}{p}
\newcommand{\Oy}{q}

\usepackage[OT2,OT1]{fontenc}
\newcommand\cyr{\fontencoding{OT2}\fontfamily{wncyr}\selectfont}

\newcommand{\tr}{{\rm tr\,}}		
\renewcommand{\Re}{{\rm Re}}
\renewcommand{\bar}{\overline}

\renewcommand\AA{{\cal A}}
\newcommand\II{{\cal I}}
\newcommand\VV{{\cal V}}
\newcommand\up{u^{\phantom\dagger}}
\newcommand\vp{v^{\phantom\dagger}}

\begin{document}


\title{\bfseries ORTHONORMAL EIGENBASES \\ {\large OVER THE} \\[1mm] OCTONIONS}

\author{
	Tevian Dray \\
	{\itshape Department of Mathematics} \\ 
	{\itshape Oregon State University} \\
	{\itshape Corvallis, OR  97331} \\
	{\ttfamily tevian{\rmfamily @}math.orst.edu} \\
		\and
	Corinne A. Manogue \\
	{\itshape Department of Physics} \\
	{\itshape Oregon State University} \\
	{\itshape Corvallis, OR  97331} \\
	{\ttfamily corinne{\rmfamily @}physics.orst.edu} \\
		\and
	Susumu Okubo
	\thanks{Supported in part by
		U. S. Department of Energy Grant No.\ DE-F402-ER 40685.}
	\\
	{\itshape Department of Physics and Astronomy} \\
	{\itshape University of Rochester} \\
	{\itshape Rochester, NY  14627} \\
	{\ttfamily okubo{\rmfamily @}pas.rochester.edu}
}

\date{\today}

\maketitle

\begin{abstract}
We previously showed that the real eigenvalues of $3\times3$ octonionic
Hermitian matrices form two separate families, each containing 3 eigenvalues,
and each leading to an orthonormal decomposition of the identity matrix, which
would normally correspond to an orthonormal basis.  We show here that it
nevertheless takes both families in order to decompose an arbitrary vector
into components, each of which is an eigenvector of the original matrix; each
vector therefore has 6 components, rather than 3.
\end{abstract}

\section{Introduction}

In previous work~\cite{Eigen,Find,OkuboPaper} we considered the real eigenvalue
problem for $3\times3$ octonionic Hermitian matrices.
\footnote{There are also eigenvalues which are not
real~\cite{NonReal,NonReal2}, whose intriguing physical properties were
discussed in~\cite{Dim,Spin}.  A related eigenvalue problem, which admits only
real eigenvalues, was discussed in~\cite{Other}.}
For this case, there are 6, rather than 3, real eigenvalues~\cite{Ogievetsky}.
We showed that these come in 2 independent families, each consisting of 3 real
eigenvalues, which satisfy a modified characteristic equation rather than the
usual one.  Furthermore, the corresponding eigenvectors are not orthogonal in
the usual sense, but do satisfy a generalized notion of orthogonality.
Finally, all such matrices admit a decomposition in terms of (the ``squares''
of) orthonormal eigenvectors, although these matrices are not idempotents
(matrices which square to themselves).

\goodbreak

It is the purpose of this paper to extend these results by showing how to
project an arbitrary 3-component octonionic vector into the eigenspaces
determined by the real eigenvalues of a given $3\times3$ octonionic Hermitian
matrix.  Somewhat surprisingly, this turns out to require all 6 eigenvalues;
generic 3-component octonionic vectors therefore have 6 vector components,
rather than 3.

\section{Octonionic Eigenvectors}

In this section, we summarize the main results of our previous
work~\cite{Eigen,Find,OkuboPaper}.

The (real) octonions $\OO$ are the noncommutative, nonassociative division
algebra generated by adjoining 7 anticommuting square roots of $-1$, with a
particular multiplication table, to the real numbers $\RR$.
\footnote{This construction can be performed over any field $\KK$, although
this will not in general yield a division algebra.  We restrict ourselves here
to the real octonions, as was done in~\cite{Eigen,Find}.  As was explicitly
verified in~\cite{OkuboPaper}, most of these results hold over an arbitrary
field, since division is only used in constructing normalized eigenvectors.}
This generalizes both the complex numbers $\CC$ and Hamilton's quaternions
$\HH$; further details can be found in~\cite{Eigen} or~\cite{GT,OkuboBook}.
We refer to a matrix as \textit{complex} or \textit{quaternionic} if its
components lie in \textit{some} complex or quaternionic subalgebra of $\OO$,
respectively.

Let $\AA$ be a $3\times3$ octonionic Hermitian matrix.  We consider solutions
of the eigenvalue equation
\begin{equation}
\AA \, v = \lambda \, v
\label{Master}
\end{equation}
where $v\in\OO^3$ and $\lambda\in\RR$.  Then $\lambda$ satisfies a modified
characteristic equation of the form
\begin{equation}
\det(\lambda I - A)
  = \lambda^3 - (\tr A) \, \lambda^2 + \sigma(A) \, \lambda - \det A
  = r
\label{Char}
\end{equation}
where the determinant is defined in terms of the Freudenthal product and where
\begin{equation}
\sigma(A) = {1\over2} \left( (\tr A)^2 - \tr (A^2) \right)
\end{equation}
The real number $r$ is a root of the quadratic equation
\begin{equation}
r^2 + 4\Phi r - | \alpha |^2 = 0
\label{Req}
\end{equation}
where both the associative 3-form $\Phi$ and the associator $\alpha$ are
totally antisymmetric functions of the non-real components of $\AA$, which are
given explicitly in Section~\ref{Tech}.  There are 3 real solutions for
$\lambda$ corresponding to each solution for $r$, which labels the families.

If the roots $r_1$ and $r_2$ of (\ref{Req}) are the same, then, over the real
octonions, $\Phi=0=\alpha$, which forces $\AA$ to be quaternionic (and
$r_1=0=r_2$).  Unless $\AA$ is complex, there are still two families of
eigenvectors, but $r$ can not be used to label them.  Rather, as described in
the next section, one family consists of the usual quaternionic eigenvectors,
while the other eigenvectors are purely octonionic, in the sense that their
components are orthogonal (as vectors in $\RR^8$) to the quaternionic
subalgebra containing the components of $\AA$.  Finally, if $\AA$ is complex,
there is only one family of eigenvectors, obtained by (right) multiplying the
usual complex eigenvectors by arbitrary octonions.

Within each family, many of the usual properties hold.  Two eigenvectors $u$,
$v$ corresponding to different eigenvalues in the same family are {\itshape
orthogonal}, but in the generalized sense that
\begin{equation}
(uu^\dagger) v = 0 = (vv^\dagger) u
\label{ortho}
\end{equation}
Even in the case of repeated eigenvalues, 3 normalized, orthogonal
eigenvectors $v_m$, with eigenvalues $\lambda_m$, can be found in each family.
Orthonormality can be expressed in the form
\begin{equation}
\II = \sum_{m=1}^3 \vp_m v_m^\dagger
\label{Idecomp}
\end{equation}
where $\II$ denotes the identity matrix, and these eigenvectors can be used to
decompose $\AA$ in the form
\begin{equation}
\AA = \sum_{m=1}^3 \lambda_m \vp_m v_m^\dagger
\label{Adecomp}
\end{equation}

\section{Quaternionic Projections}

In this section, we consider the much easier case of quaternionic matrices.
While we are primarily interested in the $3\times3$ case, all results in this
section apply to $n\times n$ quaternionic Hermitian matrices.

It is easily shown~\cite{Eigen} that the right eigenvalue problem
\begin{equation}
\AA \, v = v\,\lambda
\end{equation}
for quaternionic Hermitian matrices $\AA$ admits solutions only for
$\lambda\in\RR$.  It is straightforward to find an orthonormal basis
$\{v_n\}$ of (quaternionic) eigenvectors, which satisfies
\begin{eqnarray}
\AA v_m &=& \lambda_m v_m \\
v_m^\dagger \vp_n &=& \delta_{mn}
\end{eqnarray}
Given any quaternionic vectors $v$, $x$, the projection of $x$ along $v$ is
simply $v(v^\dagger x)$, so that any quaternionic vector has vector components
of the form
\begin{equation}
x = \sum \vp_m (v_m^\dagger x)
\end{equation}
which makes sense, since (the parentheses are not needed and) orthonormality
can be written as
\begin{equation}
\II = \sum \vp_m v_m^\dagger
\end{equation}

We now ask about the \textit{octonionic} properties of the quaternionic
matrix $\AA$.  To do this, we view the octonions through the Cayley-Dickson
process as
\begin{equation}
\OO = \HH \oplus \ell\HH
\end{equation}
where $\ell$ is any octonionic unit orthogonal (as vectors in $\RR^8$) to all
elements of $\HH$, the quaternionic subalgebra of $\OO$ containing the
components of $\AA$.
\footnote{If $\AA$ is complex, the choice of $\HH\subset\OO$ is not fixed.  In
this case, however, there is only one family of (octonionic) eigenvectors.  We
assume in what follows that $\AA$ is not complex.}

What are the \textit{octonionic} eigenvectors of $\AA$?  As discussed
in~\cite{Eigen}, if the components of $\AA$ and $v$ lie in the (same)
quaternionic subalgebra $\HH$, then we have
\begin{equation}
\AA(\ell v) = \ell (\bar\AA v)
\end{equation}
Thus, if $u$ is an eigenvector of $\bar\AA$, then $\ell u$ is an eigenvector
of $\AA$ with the same eigenvalue!  This implies that, over $\OO$, the real
eigenvalues of $\AA$ consist of its eigenvalues over $\HH$ together with the
eigenvalues of $\bar\AA$ over $\HH$.

As expected, there are two families of (octonionic) eigenvectors of $\AA$, the
usual family from the quaternionic eigenvalue problem, together with an
``extra'' family, namely $\{\ell u_m\}$ where $\{u_m\}$ are the (quaternionic)
eigenvectors of $\bar\AA$, with real eigenvalues $\mu_m$.  This latter family
also leads to a decomposition of $\AA$.  To see this, note that if $u$ is
quaternionic then
\begin{equation}
(\ell u)(\ell u)^\dagger = \ell (uu^\dagger) \bar\ell = \bar{uu^\dagger}
\end{equation}
from which it follows that for any quaternionic $u$, $x$
\begin{equation}
\left[(\ell u)(\ell u)^\dagger\right] (\ell x)
  = \bar{uu^\dagger} (\ell x)
  = \ell (uu^\dagger x)
\end{equation}
But this means that $(\ell u)(\ell u)^\dagger$ has precisely the same
projection properties over $\ell\HH$ as $uu^\dagger$ has over $\HH$.  Thus, if
$\{u_m\}$ are (quaternionic) orthonormal eigenvectors of $\bar\AA$, then
$\{\ell u_m\}$ are (octonionic) orthonormal eigenvectors of $\AA$, and in
particular
\begin{eqnarray}
\II &=& \sum (\ell \up_m) (\ell u_m^\dagger) \\
\AA &=& \sum \mu_m (\ell \up_m) (\ell u_m^\dagger)
\end{eqnarray}

Putting this all together, we can decompose an arbitrary octonionic vector
$\Vx$ in terms of the eigenvectors of $\AA$ and $\bar\AA$.  First of all,
there is a unique decomposition of $\Vx$ into a quaternionic piece and a
``purely octonionic'' piece, given by
\begin{equation}
\Vx = \Vx_1 + \ell \Vx_2
\end{equation}
where both $\Vx_1$ and $\Vx_2$ are quaternionic.  The rest is easy: Expand
$\Vx_1$ in terms of the (quaternionic) eigenvectors of $\AA$, and $\Vx_2$ in
terms of the (quaternionic) eigenvectors of $\bar\AA$.  Equivalently, expand
$\Vx_1$ in terms of $\{v_m\}$ and $\ell \Vx_2$ in terms of $\{\ell u_m\}$.
This leads to
\begin{eqnarray}
\Vx &=& \sum \vp_m (v_m^\dagger \Vx_1)
	+ \sum \left( (\ell \up_m)(\ell u_m^\dagger) \right) (\ell \Vx_2)
	\nonumber\\
    &=& \sum \vp_m v_m^\dagger \Vx_1 + \ell \sum \up_m u_m^\dagger \Vx_2
\end{eqnarray}

We have therefore succeeded in projecting $\Vx$ into vector components, each
of which is an eigenvector of $\AA$.  However, there are twice as many
components as one would have expected from the purely quaternionic problem.
Furthermore, while it is clear that $\vp_m v_m^\dagger \Vx_1$ is an
eigenvector of $\AA$ with eigenvalue $\lambda_m$, since it is a quaternionic
right multiple of $v_m$, it is rather surprising that $\ell (u_m u_m^\dagger
\Vx_1)$ is an eigenvector of $\AA$, since it is \textit{not} in general a
multiple of $\ell u_m$!

We will see in the next section that this entire structure carries over to the
fully octonionic setting, but that the identification of the two families is
not as easy.

\section{Decompositions of the Octonions}
\label{Tech}

In this section, we summarize and interpret some of the more technical results
from \cite{OkuboPaper}, which we will need to prove our main theorems
in Section~\ref{Proofs}.  These results are interesting in their own right,
however, as they lead to the consideration of several intriguing ways of
decomposing the octonions into subspaces.

Write the components of the octonionic Hermitian matrix $\AA$ as
\begin{equation}
\AA = \pmatrix{\Ap& \Aa& \bar{\Ab}\cr
	\bar{\Aa}& \Am& \Ac\cr \Ab& \bar{\Ac}& \An\cr}
\end{equation}
with $\Ap,\Am,\An\in\RR$ and $\Aa,\Ab,\Ac\in\OO$.  With these conventions,
the associative 3-form \cite{Harvey} is given by
\begin{equation}
\Phi
  = \Re (\Aa \times \Ab \times \Ac)
  = {1\over2} \Re \big( \Aa(\bar{\Ab}\Ac) - \Ac(\bar{\Ab}\Aa) \big)
\end{equation}
and the associator is given by
\begin{equation}
\alpha =
[\Aa,\Ab,\Ac] = (\Aa\Ab)\Ac-\Aa(\Ab\Ac)
\end{equation}
Viewing the octonions as the real vector space $\RR^8$, we introduce the
subspace $\TT\subset\OO$ spanned by the elements of $\AA$.  In other words,
$\TT$ is the real vector space consisting of all (real) linear combinations of
the components of $\AA$, that is
\begin{equation}
\TT = \langle 1,\Aa,\Ab,\Ac \rangle
\end{equation}
Since the quaternionic case was covered in the previous section, we will
assume throughout the remainder of the paper that $\alpha\ne0$, which in turn
implies that $\TT$ has dimension 4 (as a real vector space).

Define the \textit{characteristic operator} $K$ associated with $\AA$ by
\begin{equation}
K[\Vx] = \AA \big( \AA (\AA\Vx) \big) - (\tr\AA) \big( \AA (\AA\Vx) \big)
		+ \sigma(\AA) (\AA\Vx) - (\det\AA) \Vx
\label{Keq}
\end{equation}
where $\Vx\in\OO^3$.  As we may readily verify, the operator $K$ is diagonal,
and can thus be reinterpreted as an operator acting on single octonions.  This
operator, also called $K$, can be written in the form
\begin{equation}
K[\Ox]
  = \Ac(\Ab(\Aa\Ox)) + \bar{\Aa}(\bar{\Ab}(\bar{\Ac}\Ox)) - 
    \Big( \Ac(\Ab\Aa) + \bar{\Aa}(\bar{\Ab}\bar{\Ac}) \Big) \Ox
\end{equation}
for any $\Ox\in\OO$.  In this form, $K$ is the same as the operator introduced
in~\cite{OkuboPaper}, where it was shown to satisfy the same quadratic
equation as~$r$, that is (compare (\ref{Req}))
\begin{equation}
K^2 + 4\Phi K - | \alpha |^2 = 0
\label{Kid}
\end{equation}
In fact, comparing (\ref{Keq}) with (\ref{Char}), we see that if $v$ is an
eigenvector of $\AA$, so that (\ref{Master}) holds, then
\begin{equation}
K[v] = rv
\end{equation}
with $r$ a solution of (\ref{Req}).  Thus, the operator $K$ can be used to
distinguish the 2 families of eigenvectors, labeled by the 2 possible values
of $r$.
\footnote{Recall that, over the real octonions, the two solutions of
(\ref{Req}) are distinct unless $\Phi=0=\alpha$, in which case $\AA$ is
quaternionic and can be handled by the methods of the previous section, and
$K$ is not needed to distinguish the two families.}

Using results from~\cite{OkuboPaper}, we can show that $K$ takes a
particularly simple form on $\TT$, namely
\begin{equation}
K[t] = t\alpha \qquad (t\in\TT)
\end{equation}
Since
\begin{equation}
\TT^\perp \equiv \TT\alpha
\end{equation}
it follows using~(\ref{Kid}) that
\begin{equation}
K[u] = -u (\alpha+4\Phi) \qquad (u\in\TT^\perp)
\end{equation}
An intriguing property of the decomposition
\begin{equation}
\OO = \TT \oplus \TT\alpha
\end{equation}
is that it is almost of the Cayley-Dickson form.  Direct computation
establishes the multiplication table (for $t_1,t_2\in\TT$)
\begin{eqnarray}
t_1(t_2\alpha) &=& (t_2t_1)\alpha \\
(t_1\alpha)t_2 &=& (t_1\bar{t_2})\alpha \\
(t_1\alpha)(t_2\alpha) &=& -\bar{t_2}t_1|\alpha|^2
\end{eqnarray}
which is of the Cayley-Dickson form.  However, this is not a true
Cayley-Dickson extension of $\TT$, since $\TT$ is not closed under
multiplication.

Using these ideas, it can be shown (Proposition 3.1 of~\cite{OkuboPaper}) that
the eigenspace of $K$ with eigenvalue $r$ is precisely $\TT(r+4\Phi+\alpha)$,
that is
\begin{equation}
K[q] = r_m q \Longleftrightarrow \exists t\in\TT : q = t(r_m+4\Phi+\alpha)
\label{Prop31}
\end{equation}
where the solutions of (\ref{Req}) are denoted $r_m$, with $m=1,2$.  If we
define
\footnote{Since $r_1+r_2+4\Phi=0$, the denominator can only vanish if
$r_1=r_2$, or equivalently $\tau=0$ in the notation of~\cite{OkuboPaper};
this cannot happen here since we are assuming $\alpha\ne0$.}
\begin{equation}
s_m = \frac{r_m+4\Phi+\alpha}{2(r_m+2\Phi)}
\end{equation}
then the decomposition of $\OO$ into its eigenspaces under $K$ takes the form
\begin{equation}
\OO = \TT s_1 \oplus \TT s_2
\end{equation}
and the normalization is chosen such that
\begin{equation}
s_1 + s_2 = 1
\end{equation}
The eigenspaces $\TT_m:=\TT s_m$ are orthogonal, since $K$ is self-adjoint.
Explicitly, the inner product in $\RR^8$ takes the form
\begin{equation}
2 \Ox\cdot \Oy
  = \Ox\bar{\Oy} + \Oy\bar{\Ox}
  = \bar{\Ox}\Oy + \bar{\Oy}\Ox
  \qquad (\Ox,\Oy\in\OO)
\end{equation}
and $K$ satisfies
\begin{equation}
K[\Ox] \cdot \Oy = \Ox \cdot K[\Oy]
\end{equation}
It is now straightforward to define projection operators
\begin{equation}
K_m = {K+r_m+4\Phi \over 2(r_m+2\Phi)}
\end{equation}
from $\OO$ to $\TT_m$, which satisfy
\begin{equation}
K_1 + K_2 = 1 \qquad K_1 K_2 = 0 = K_2 K_1
  \qquad K_1^2 = K_1 \qquad K_2^2 = K_2
\end{equation}
Finally, since $\bar{s_1} s_2$ is proportional to $\alpha$, we also have
\begin{equation}
\TT_2 \equiv \TT_1 \alpha
\end{equation}

The following two results are key to understanding the family structure.
First, assertions (v) and (vi) of Proposition 3.2 of~\cite{OkuboPaper} imply
\begin{equation}
\Ox,\Oy\in\TT_m \quad\Longrightarrow\quad \Ox\bar{\Oy}\in\TT
\label{Prop32}
\end{equation}
Second, as we also see from Proposition 4.1 of~\cite{OkuboPaper}, for any
\hbox{$\Ox_1,\Ox_2\in\TT$} and $\Oy\in\TT_m$, we have
\begin{equation}
[\Ox_1,\Ox_2,\Oy] = \Ox_3\Oy
\label{Prop41}
\end{equation}
where $\Ox_3\in\OO$ depends on $\Ox_1,\Ox_2$ (and $m$), but does \textit{not}
depend on $\Oy$.  These two results are the primary ingredients in the proofs
in Section~\ref{Proofs} of our main results.

\section{Octonionic Projections}

We now return to the discussion of octonionic eigenvectors, showing how to
project any vector into the 6 eigenspaces they determine.  We begin with the
following result.

\goodbreak
\textit{{\bfseries Theorem 1:}
Let a $3\times3$ octonionic Hermitian matrix $\AA$ be given, and suppose that
$v\in\OO^3$ is an eigenvector of $\AA$, that is, $\AA v=\lambda v$ with
$\lambda\in\RR$.  Suppose further that $K[v]=rv$, with $K$ as in~(\ref{Keq}),
and let $\Vy\in\OO^3$ be any other vector satisfying $K[\Vy]=r\Vy$ (for the
same $r$).
Then}
\begin{equation}
(vv^\dagger) \left( (vv^\dagger)\Vy \right) = 
 (v^\dagger v) \left( (vv^\dagger)\Vy \right)
\label{result}
\end{equation}

\textit{{\bfseries Proof:}
See Section~\ref{Proofs}.}

In other words, $(vv^\dagger)\Vy$ is an eigenvector of a new $3\times3$
octonionic Hermitian matrix, namely $vv^\dagger$; if $v$ is normalized, the
eigenvalue is~$1$.  Note, however, that $(vv^\dagger) \Vy$ is not necessarily
a multiple of $v$!  If $u$ is another eigenvector of $\AA$ in the same family,
so that $K[u]=ru$ and $\AA u=\mu u$ and where we assume $\lambda\ne\mu\in\RR$,
then $u$ and $v$ are orthogonal (in the sense~(\ref{ortho})).  This in turn
means that $u$ and $(vv^\dagger)\Vy$ are eigenvectors of the matrix
$vv^\dagger$ with different eigenvalues, and hence must themselves be
orthogonal, that is
\begin{equation}
uu^\dagger \left( (vv^\dagger)\Vy \right) = 0
\end{equation}
This shows that, {\it restricted to the appropriate $K$ eigenspace},
$uu^\dagger$ and $vv^\dagger$ are projection operators!  Extending this to a
family $\{u,v,w\}$ of orthonormal eigenvectors, we obtain using~(\ref{Adecomp})
\begin{equation}
\AA \left( (vv^\dagger)\Vy \right) = \lambda \left( (vv^\dagger)\Vy \right)
\end{equation}
Furthermore, (\ref{Idecomp}) leads to
\begin{equation}
\Vy = (uu^\dagger)\Vy + (vv^\dagger)\Vy + (ww^\dagger)\Vy
\label{ydecomp}
\end{equation}
which expresses $\Vy$ as the sum of its projections along the eigenvectors.

Given an arbitrary vector $\Vx\in\OO^3$, we can use the projection operators
$K_m$ to project $\Vx$ into the eigenspaces of $K$, by defining
\begin{equation}
\Vx_n = K_n[\Vx]
\end{equation}
and noting that
\begin{equation}
\Vx = \Vx_1 + \Vx_2
\end{equation}
We can then use (\ref{ydecomp}) to expand each of $\Vx_m$ along the
appropriate eigenvectors.  Putting this all together, any vector $\Vx$ can be
expanded into \textit{six} terms, consisting first of the projections into the
2 eigenspaces of $K$, then decomposing each of these separately into 3
eigenspaces of~$\AA$.

\section{Families of Eigenvectors}

Is the family structure a property of the matrix $\AA$, or of its
eigenvectors?  We show in this section that it is the latter.

Decompose $\AA$ in terms of an orthonormal eigenbasis as in (\ref{Adecomp}),
yielding
\begin{equation}
\AA=\lambda uu^\dagger+\mu vv^\dagger+\nu ww^\dagger
\label{Matrix}
\end{equation}
Then, by the results of the last section, all of $u$, $v$, $w$ have the same
eigenvalue under $K$, namely
\begin{equation}
r = \lambda\mu\nu - \det\AA
\end{equation}
Note that this result holds regardless of the values of $\lambda$, $\mu$,
$\nu$; it is a property of the orthonormal ``basis'' $\{u,v,w\}$, not of a
particular matrix $\AA$.

Consider the special case where $\AA=vv^\dagger$, with $v^\dagger v=1$.
Then $v$ is an eigenvector of $\AA$ with eigenvalue $1$, and
\begin{equation}
K[v] = -(\det\AA) v
\end{equation}
since
\begin{equation}
\tr\AA = 1
\qquad
\sigma(\AA) = 0
\end{equation}
In this case, (\ref{result}) simplifies to
\begin{equation}
K[\Vy] = -(\det vv^\dagger) \Vy \Longleftrightarrow
vv^\dagger \left( (vv^\dagger)\Vy \right) = (vv^\dagger)\Vy
\end{equation}
since the converse follows immediately from the definition of $K$ (with
$\AA=vv^\dagger$).  Special cases of vectors $u$ satisfying this condition are
$u$ orthogonal to $v$, so that
\begin{equation}
\AA u = (vv^\dagger)u = 0
\end{equation}
or
\begin{equation}
uu^\dagger = s \, vv^\dagger
\label{Phase}
\end{equation}
with $s\in\RR$, in which case
\begin{equation}
\AA u
  = (vv^\dagger) u
  = \frac{1}{s} (uu^\dagger) u
  = \frac{1}{s} u (u^\dagger u)
  = u
\end{equation}
The decomposition (\ref{ydecomp}) shows that any $\Vy$
satisfying
\begin{equation}
K[\Vy] = -(\det\AA) \Vy
\label{yeq}
\end{equation}
can be written as the sum of terms of one of these two types.

If $vv^\dagger$ is complex, then $v$ determines a single ``family'' containing
all vectors.  Furthermore, $v$ must then be a (right) octonionic multiple of a
complex vector $v_0$; the ``family'' structure in this case corresponds to
completing $v_0$ as usual to a complex orthonormal basis.  However, there are
also other families containing $v$.  (If $v$ is a \textit{real} eigenvector of
$\AA$, then $\AA$ must be real or complex, and there is then no nontrivial
family containing $v$.)

If $vv^\dagger$ is not complex, there is an $8$-parameter family of vectors
orthogonal to $v$, and a $4$-parameter family of vectors satisfying
(\ref{Phase}).  We therefore see that $v$ alone determines a 12-parameter
family of vectors satisfying both (\ref{result}) and (\ref{yeq}).  This family
consists precisely of the (real) linear combinations of a family of
eigenvectors of any matrix of the form (\ref{Matrix}).  In other words, a
family of eigenvectors is just an orthonormal ``basis'', but each such basis
spans only half of $\OO^3$.

We can complete this picture by showing that each vector $v$ (with
$vv^\dagger$ not complex) lies in exactly one such family.  To do this, we
first remove the requirement that $v$ be an eigenvector from (\ref{result}).
Specifically, we have

\goodbreak
\textit{{\bfseries Theorem 2:}
Let $K$ be the operator constructed via (\ref{Keq}) from a given $3\times3$
octonionic Hermitian matrix $\AA$, and let $\Vy,\Vz\in\OO^3$.  Suppose
$K[\Vy]=r\Vy$ and $K[\Vz]=r\Vz$ with $r\in\RR$.  Then}
\begin{equation}
(\Vy\Vy^\dagger) \left( (\Vy\Vy^\dagger)\Vz \right) = 
 (\Vy^\dagger \Vy) \left( (\Vy\Vy^\dagger)\Vz \right)
\label{zeq}
\end{equation}

\textit{{\bfseries Proof:}
See Section~\ref{Proofs}.}

This finally allows us to argue that if $\Vy$ is in the family defined by $v$,
by virtue of satisfying either (and hence both) of (\ref{result}) and
(\ref{yeq}), then $v$ is also in the family defined by $\Vy$.  This follows
immediately from Theorem 2, since (\ref{yeq}) is enough to satisfy the
condition of the theorem, and (\ref{zeq}) (with $\Vz=v$) ensures that $v$ is
in the family defined by $\Vy$.

More generally, Theorem 2 shows that $u$, $w$ belong to the same family (with
$u$ normalized by $u^\dagger u=1$) if and only if
\begin{equation}
uu^\dagger \left((uu^\dagger)w\right) = \left((uu^\dagger)w\right)
\label{weq}
\end{equation}
For, (\ref{weq}) shows that $w$ is in the family defined by $u$.  If $u$ is in
the family defined by $v$, say, then the previous argument shows that $v$ is
also in the family defined by $u$.  But then $v$ and $w$ are both in the
family defined by $u$, and the theorem shows that they must each be in the
family defined by the other.  This shows that any family which contains $u$
must also contain~$w$.  The converse, namely that $u$, $w$ in the same family
implies (\ref{weq}), follows immediately from the theorem.

\section{Discussion}
\label{Discussion}

As made clear in Section~\ref{Tech}, the family structure, which was originally
discovered on $\OO^3$ arising from the eigenvectors of $\AA$, can be viewed as
a property of the octonions themselves.  From this point of view, the only
purpose of $\AA$ is in determining the octonions $\Aa,\Ab,\Ac$ used to define
$\TT$.

So suppose $\Aa,\Ab,\Ac\in\OO$ are given, and further assume that
$\alpha=[\Aa,\Ab,\Ac]\ne0$.  Does this decomposition of $\OO$ depend on the
choice of $\Aa$, $\Ab$, $\Ac$, or merely on $\TT$?  Suppose $A,B,C\in\TT$, so
that
\begin{equation}
A = A_0 + A_1 \Aa + A_2 \Ab + A_3 \Ac
\end{equation}
with $A_m\in\RR$, and similarly for $B$ and $C$.  If we write $\vec{A}$ for
$(A_1,A_2,A_3)\in\RR^3$, etc., then
\begin{eqnarray}
[A,B,C] &=& (\vec{A},\vec{B},\vec{C}) \> [\Aa,\Ab,\Ac] \\
\Phi(A,B,C) &=& (\vec{A},\vec{B},\vec{C}) \> \Phi(\Aa,\Ab,\Ac)
\end{eqnarray}
where $(\vec{A},\vec{B},\vec{C}) = (\vec{A}\times\vec{B})\cdot\vec{C}$ is the
vector triple product, assumed to be nonzero.  Since both $\alpha$ and $\Phi$
scale with the same factor, so does $r$ --- which means that the $s_m$ are
unchanged!  The decomposition of $\OO$ thus depends only on the subspace
$\TT$.

In terms of matrices, this means that any matrix with non-real components
$A,B,C\in\TT$ satisfying $(\vec{A},\vec{B},\vec{C})\ne0$ determines the same
family structure as the original matrix $\AA$.  The 2 families determined by
$\AA$ are just
\begin{equation}
\OO^3 = \FF_1 \oplus \FF_2 := (\TT_1)^3 \oplus (\TT_2)^3
\end{equation}
An essential ingredient in this picture is the result~(\ref{Prop32}), which
ensures that, for {\it any\/} vector $v\in\FF_m$, the components of
$vv^\dagger$ are in $\TT$!

We conclude by emphasizing that the subspaces $\TT_m$ have many of the
properties of the quaternions: They are 4-dimensional, they lead to an almost
Cayley-Dickson multiplication table for the octonions, and, our main result,
they have just enough associativity to enable arbitrary vectors to be
projected into eigenspaces.

\newpage
\goodbreak
\section{Proofs of the Theorems}
\label{Proofs}

We give here the proofs of our two theorems.  This section can be skipped on
first reading; all results have also been independently derived using {\sl
Mathematica}.

\goodbreak
\textit{{\bfseries Theorem 1:}
Let a $3\times3$ octonionic Hermitian matrix $\AA$ be given, and suppose that
$v\in\OO^3$ is an eigenvector of $\AA$, that is, $\AA v=\lambda v$ with
$\lambda\in\RR$.  Suppose further that $K[v]=rv$, with $K$ as in~(\ref{Keq}),
and let $\Vy\in\OO^3$ be any other vector satisfying $K[\Vy]=r\Vy$ (for the
same $r$).
Then}
\begin{eqnarray*}
(vv^\dagger) \left( (vv^\dagger)\Vy \right) = 
 (v^\dagger v) \left( (vv^\dagger)\Vy \right)
\end{eqnarray*}

\textit{{\bfseries Proof:}}
By~(\ref{Prop31}), $v,y\in(\TT_m)^3$, and by~(\ref{Prop32}) the components of
$vv^\dagger$ are in $\TT$.  Thus, introducing a natural matrix associator,
(\ref{Prop41}) implies that
\begin{equation}
(vv^\dagger) \left( (vv^\dagger)\Vy \right)
  = \left( (vv^\dagger) (vv^\dagger) \right) \Vy
	+ [vv^\dagger,vv^\dagger,\Vy]
  = \VV\Vy
\end{equation}
where $\VV$ does {\itshape not} depend on $\Vy$!  But alternativity implies
that the vector associator satisfies~\cite{NonReal}
\begin{equation}
[v,v^\dagger,v] \equiv 0
\end{equation}
which in turn implies that
\begin{equation}
(vv^\dagger) \left( (vv^\dagger)v \right)
  = (v^\dagger v) \left( (vv^\dagger)v \right)
\label{Vid}
\end{equation}
(and the right-hand-side simplifies even further).  We want to show that
\begin{equation}
\VV = (v^\dagger v) (vv^\dagger)
\label{Veq}
\end{equation}
but~(\ref{Vid}) in fact only shows that
\begin{equation}
\VV = (v^\dagger v) (vv^\dagger) + \VV_1; \qquad \VV_1 v = 0
\end{equation}
However, if $u$ is a different eigenvector of $\AA$ in the same family, that
is, if
\begin{eqnarray}
\AA u &=& \mu u \\
 K[u] &=& r_m u
\end{eqnarray}
where we assume without loss of generality that $\mu\ne\lambda$, then the
orthogonality~(\ref{ortho}) of $u$ and $v$ implies that
\begin{equation}
(vv^\dagger) \left( (vv^\dagger)u \right)
  = (v^\dagger v) \left( (vv^\dagger)u \right)
\end{equation}
(since each side is $0$).  Since the right-hand-side is $\VV u$, we must also
have
\begin{equation}
\VV_1 u = 0
\end{equation}
so that $\VV_1=0$ when acting on any eigenvector of $\AA$ (in the same family
as $v$).  But this means that $\AA$ and $\AA+\VV_1$ have the same
decomposition and are therefore equal, which forces $\VV_1=0$, and
establishes~(\ref{Veq}).
\null\hfill$\blacksquare$

\goodbreak
\textit{{\bfseries Theorem 2:}
Let $K$ be the operator constructed via (\ref{Keq}) from a given $3\times3$
octonionic Hermitian matrix $\AA$, and let $\Vy,\Vz\in\OO^3$.  Suppose
$K[\Vy]=r\Vy$ and $K[\Vz]=r\Vz$ with $r\in\RR$.  Then}
\begin{eqnarray*}
(\Vy\Vy^\dagger) \left( (\Vy\Vy^\dagger)\Vz \right) = 
 (\Vy^\dagger \Vy) \left( (\Vy\Vy^\dagger)\Vz \right)
\end{eqnarray*}

\textit{{\bfseries Proof:}}
As in the proof of the Theorem~1, $\Vy,\Vz\in(\TT_m)^3$
by~(\ref{Prop31}), and the components of $\Vy\Vy^\dagger$ are in $\TT$.  We
can therefore write
\begin{equation}
\Vy\Vy^\dagger
  = \pmatrix{d_1& t_3& \bar{t_2}\cr
	\bar{t_3}& d_2& t_1\cr t_2& \bar{t_1}& d_3\cr}
\label{Ycomps}
\end{equation}
with $d_n\in\RR$ and $t_n\in\OO$.  Introducing the components
$\Vy_n,\Vz_n\in\TT_m$ of $\Vy,\Vz$, direct computation using alternativity
shows that
\begin{equation}
t_2 (t_3 \Vy_2)
  = \Vy_3 |\Vy_1|^2 |\Vy_2|^2
  = d_1 \bar{t_1} \Vy_2
\end{equation}
and cyclic permutations, where we have used the fact that
\begin{eqnarray}
d_1 &=& |t_3|^2 \label{outer1}\\
t_3 &=& \Vy_1 \bar{\Vy_2}\label{outer2}
\end{eqnarray}
and so forth.  But (\ref{Prop41}) implies that
\begin{equation}
t_2 (t_3 q)
  = (t_2 t_3) q + [t_2,t_3,q]
  = p q
\end{equation}
for some $p\in\OO$ which is independent of $q\in\TT_m$, and we have therefore
shown that
\begin{equation}
p = d_1 \bar{t_1}
\end{equation}
In particular, we must have
\begin{equation}
t_2 (t_3 \Vz_2) = d_1 \bar{t_1} \Vz_2
\label{assoc1}
\end{equation}
and cyclic permutations.  A similar argument establishes
\begin{equation}
\bar{t_1} (\bar{t_3} q) = d_2 t_2 q
\label{assoc2}
\end{equation}
for any $q\in\TT_m$, together with its cyclic permutations.

We are now ready to compute both sides of~(\ref{zeq}) explicitly.  We have
first of all
\begin{equation}
(\Vy\Vy^\dagger) \Vz
  = \pmatrix{d_1\Vz_1+t_3\Vz_2+\bar{t_2}\Vz_3\cr
		\bar{t_3}\Vz_1+d_2\Vz_2+t_1\Vz_3\cr
		t_2\Vz_1+\bar{t_1}\Vz_2+d_3\Vz_3\cr}
\end{equation}
Multiply on the left by $\Vy\Vy^\dagger$, and consider for simplicity only the
last component, which is
\begin{eqnarray}
\lefteqn{
t_2 (d_1\Vz_1+t_3\Vz_2+\bar{t_2}\Vz_3)
  + \bar{t_1} (\bar{t_3}\Vz_1+d_2\Vz_2+t_1\Vz_3)
  + d_3 (t_2\Vz_1+\bar{t_1}\Vz_2+d_3\Vz_3)
}\nonumber\\
  &=& (d_1+d_3) t_2 \Vz_1 + \bar{t_1}(\bar{t_3}\Vz_1)
	+ t_2(t_3\Vz_2) + (d_2+d_3) \bar{t_1} \Vz_2
	+ (|t_1|^2+|t_2|^2+d_3^2) \Vz_3 \\
  &=& (d_1+d_2+d_3) (t_2\Vz_1+\bar{t_1}\Vz_2+d_3\Vz_3) \nonumber
\end{eqnarray}
where we have used (\ref{outer1}), (\ref{outer2}), (\ref{assoc1}), and
(\ref{assoc2}) in the last equality.  Repeating this for the remaining
components, and noticing that
\begin{equation}
d_1 + d_2 + d_3 = \tr(\Vy\Vy^\dagger) = \Vy^\dagger\Vy
\end{equation}
completes the proof.
\null\hfill$\blacksquare$

Theorem~1 clearly follows from Theorem~2.  But one can also argue that the
reverse is true, thus avoiding the detailed computation given above.  The
argument goes as follows: $\Vy\Vy^\dagger$ is itself a Hermitian matrix with
components in $\TT$, and, by the argument given in Section~\ref{Discussion},
it determines the same family structure as the given matrix $\AA$, and can
therefore be used instead of $\AA$ in applying Theorem~1.

\newpage


\end{document}